# Contradictions and narrowness of views in "The fables of Ishango, or the irresistible temptation of mathematical fiction", answers and updates

Vladimir Pletser ([1]), Dirk Huylebrouck ([2])

*(1) Sciences Department, Manned Spaceflight and Operations Directorate, European Space Research and Technology Centre, European Space Agency, Noordwijk, The Netherlands;*
*Present address: Technology and Engineering Centre for Space Utilization,*
*Chinese Academy of Sciences, Beijing, China*
*(2) Faculty of Architecture, University of Leuven, KUL Brussels-Gent, Belgium*

Abstract

We answer to criticisms of O. Keller about our interpretation work on the Ishango rod, the oldest mathematical tool of humankind. Our hypothesis, that is widely accepted, is that this prehistoric rod is the first mankind manifestation of a basic arithmetic intention, with simple arithmetic operations and possibly showing passages between 10 and 12 bases.

―――――――

# Contradictions et étroitesse de vues dans « Fables d'Ishango, ou l'irrésistible tentation de la mathématique-fiction », réponses et mises au point

Vladimir Pletser([1]), Dirk Huylebrouck([2])

*(1) Département Sciences, Directorat Vols spatiaux habités et Opérations, Centre de Recherches et de Technologies spatiales, Agence spatiale européenne, Noordwijk, Pays-Bas ;*
*Actuellement au Centre de Technologie et d'Ingénierie pour l'Utilisation Spatiale de l'Académie chinoise des Sciences, Pékin, Chine*
*(2) Faculty for Architecture, University of Leuven, KUL Brussels-Gent, Belgium*

Abstract

Nous répondons aux critiques d'O. Keller au sujet de notre travail d'interprétation du bâton d'Ishango, l'outil mathématique le plus ancien de l'humanité. Notre hypothèse, qui est largement acceptée, est que ce bâton préhistorique est la première manifestation de l'humanité d'une intention arithmétique de base, avec des opérations arithmétiques simples et utilisant probablement des passages entre les bases 10 et 12.

# Contradictions et étroitesse de vues dans « Fables d'Ishango, ou l'irrésistible tentation de la mathématique-fiction », réponses et mises au point

Vladimir Pletser($^1$), Dirk Huylebrouck($^2$)

*(1) Département Sciences, Directorat Vols spatiaux habités et Opérations, Centre de Recherches et de Technologies spatiales, Agence spatiale européenne, Noordwijk, Pays-Bas ;
Actuellement au Centre de Technologie et d'Ingénierie pour l'Utilisation Spatiale de l'Académie chinoise des Sciences, Pékin, Chine
(2) Faculty for Architecture, University of Leuven, KUL Brussels-Gent, Belgium*

**Préliminaire : Réponse à une rumeur : de quoi s'agit-il ?**

*Une découverte faite au milieu du 20e siècle à l'Est de la République Démocratique du Congo, au cœur de la forêt des Virunga, forêt déclarée patrimoine de l'humanité par l'Unesco, d'un os comportant des encoches daté d'il y a plus de 20.000 ans balise encore notre savoir et se mue en symbole.*

*Les auteurs des lignes qui suivent, à l'instar de beaucoup de spécialistes aux qualifications certaines, ont émis leurs hypothèses. Le dictionnaire « Le Robert » définit le mot hypothèse dans un sens mathématique : « proposition admise comme donnée d'un problème » et en langage courant : « conjecture concernant l'explication ou la possibilité d'un événement ». Il précise que le sens opposé équivaut à : « conclusion, certitude, évidence ».*

*Chaque découverte bénéficie d'un acte de naissance constitué par des publications spécialisées donnant lieu à des congrès scientifiques internationaux où des experts patentés donnent leur avis. C'est ainsi qu'évolue et s'évalue son impact.*

*Mais voilà que parallèlement à ces avis argumentés circule la rumeur dénigrante. Très à l'aise dans l'utilisation de celle-ci et servi par un vocabulaire relâché, Mr Olivier Keller se répandant dans un blog, acquiert ainsi une certaine notoriété en discréditant la pensée de scientifiques plus que nuancés et prudents par profession.*

*Fallait-il l'ignorer parce que, toujours absent dans les cénacles et publications appropriées, aucune réponse immédiate et aucun débat direct avec l'auteur de ce blog n'étaient possibles ?*

*Nous avons demandé aux deux auteurs incriminés de répondre à ces critiques. Ils ont choisi de suivre à la trace chacune des 28 allégations de ce blog, d'y répondre point par point.*

*Gabriel Castel, août 2015
Comité Ishango Milele
www.ishango-milele.com*

**Introduction :**

**Les hypothèses**

Depuis la fin des années 1990, les auteurs ont présenté et expliqué les hypothèses suivantes concernant le bâton d'Ishango, découvert par l'archéologue Jean de Heinzelin à la fin des années 1950 :

1) Le bâton d'Ishango, daté de 22.000 ans, peut être considéré comme le plus ancien outil mathématique de l'humanité car l'arrangement des encoches sur trois colonnes suggère une intention arithmétique.

2) De plus, il apparaît que plusieurs bases sont utilisées dans cette arithmétique élémentaire : la base 10 et la base 12 avec ses sous-multiples 3, 4 et 6. L'arrangement géométrique des encoches dans les différents groupements sur les trois colonnes permet de faire d'autres opérations arithmétiques élémentaires.

Ces hypothèses ont été présentées dans plusieurs publications [1-18] et dans des congrès internationaux sur les mathématiques et l'ethno-mathématique.

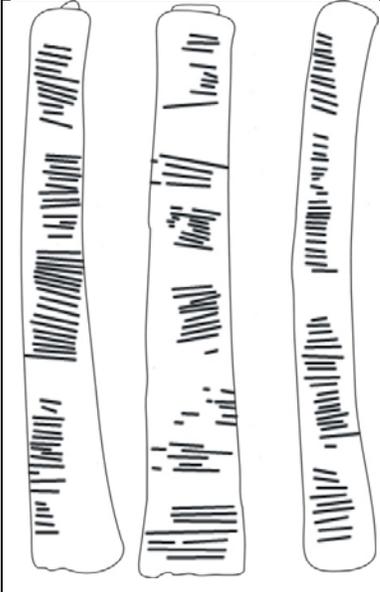

*Schéma du premier bâton d'Ishango (basé sur un dessin Wikimedia Commons) et les nombres d'encoches correspondant des trois colonnes.*

Les auteurs ont exploré en détail les relations arithmétiques entre les nombres d'encoches. On sait par exemple que les sommes des nombres dans les colonnes de gauche à droite sont 60, 48 et 60, que les nombres de la colonne de gauche peuvent être vus comme des nombres premiers, que des duplications sont vues dans la colonne du milieu, et que les nombres de la colonne de droite peuvent être interprétés comme 10 et 20 plus ou moins 1.

**La conclusion.**

Les auteurs ont proposé diverses hypothèses, et rejeté d'autres (comme celle des nombres premiers), arrivant à une conclusion peu spectaculaire : il s'agit probablement d'un outil qui dénombre des événements ou des choses, notés par quelqu'un qui mélangeait les bases 10 et 12.

Toutefois, les auteurs n'ont pas voulu entrer dans les discussions quant à l'utilisation de ce bâton : est-ce un jeu arithmétique, comme le pensait l'archéologue Jean de Heinzelin, son découvreur ; ou un calendrier, suivant l'archéologue Marshack ; ou un outil de comptage tout simplement ; ou autre chose encore ? Les auteurs ont laissé cette discussion aux ethnologues et anthropologues.

**Une critique peu académique**

Olivier Keller, dans son unique publication [19], a critiqué, voire ridiculisé, le travail des auteurs concernant ce bâton. Il leur a fallu quelque temps avant qu'ils décident de répondre à ces critiques. En effet, les travaux de ce critique n'ont jamais été publiés dans des journaux internationaux, ni dans des revues 'évaluées par des pairs' (*peer reviewed*). De plus, ce critique n'a jamais participé à aucun congrès scientifique international, où il aurait pu défendre son point de vue, ce qui fait partie de la démarche scientifique.

De nos jours, il est impossible de répondre à toutes les critiques sur des blogs ou sites internet florissant ici et là parmi certains groupes de pression ou défendant certaines convictions fanatiques ou politiques. Sans vouloir polémiquer, les auteurs veulent présenter les quelques arguments suivants pour réfuter une à une ces critiques, en italique ci-dessous et reproduites telles quelles de [19].

**<u>Les réponses :</u>**

**L'os, l'outil et son âge**

1. *« … le premier des deux os d'Ishango, est parvenu à la célébrité en se faisant passer pour un texte scientifique marquant… »*

Étrange formulation : un os se faisant passer pour un texte scientifique…

2. *« Un fragment de quartz encastré à l'une des extrémités montre qu'il s'agit d'un manche d'outil ; … »*

Cette hypothèse est soutenue par les archéologues de renom ayant étudié le bâton : Jean de Heinzelin (B), Alison Brooks (USA), John Yellen (USA), Els Cornelissen (B), et est également mentionnée sur le site de l'Institut des Sciences Naturelles de Bruxelles où l'objet est conservé.

3. *« …on le date ordinairement de 20.000 ans avant nos jours »*

La datation admise est de 22.000 ans. Cette datation est obtenue par la méthode du carbone 14, et confirmée par plusieurs autres méthodes archéologiques, comme expliqué par Alison Brooks. La civilisation d'Ishango pourrait même être âgée de 90.000 ans, suivant les travaux de J. Yellen [20].

**La suspicion sur les encoches**

4. *« … plusieurs [encoches] sont effacées ou à peine visibles, ce qui rend déjà suspecte a priori toute interprétation fondée sur leur dénombrement. »*

Il conviendrait de conseiller à ce critique de ne pas se fonder sur l'observation de photographies ou de reproductions, mais de consulter l'ouvrage original de Jean de Heinzelin [21, pp 64-70] décrivant les circonstances de sa découverte, à moins bien sûr qu'il ne faille aussi soupçonner l'illustre archéologue de malhonnêteté intellectuelle. De plus, le technicien et collaborateur de Heinzelin, Marcel Spinglaer, fut un spécialiste reconnu qui ne peut certainement pas être soupçonné.

5. « *Prenons la colonne du milieu : d'après l'auteur, 3 serait doublé en 6, 4 en 8 et 5 en 10. Mais le 5 et le 10 sont douteux : l'un des paquets de 5 est sérieusement illisible, et le 10 pourrait être en réalité 9. De plus, on ne comprend pas pourquoi, en cas de duplication de 5, 3 et 4, le groupe de cinq encoches serait figuré deux fois, contrairement au groupe de trois et de quatre qui ne le sont qu'une seule fois.* »

Où est le problème ? Ne pourrait-on y voir tout simplement deux opérations de duplication (3 => 6 et 4 => 8) et une opération d'addition (10 = 5+5) ?

6. « *Et quel est le rôle de ce 7, qui n'est ni duplication ni doublé ? À moins qu'il ne faille lire au bas de la colonne du milieu non pas 10, 5, 5 et 7, mais 10, 4, 5 et 7, auquel cas nous aurions le double de 7 avec 10+4 et le double de 5 avec 10.* »

Est-ce que l'auteur de ces lignes ne se livre pas lui-même au jeu de supposition qu'il décrie si fort par ailleurs ? Néanmoins, l'interprétation proposée est intéressante et nous en laissons la paternité et la responsabilité à son auteur.

7. « *Il est donc clair qu'à partir du moment où l'on a décidé que les paquets d'encoches sont des nombres, …* »

Il est exact que l'hypothèse de départ de toute interprétation mathématique de ce bâton d'Ishango est l'association à un groupe d'encoches du nombre d'encoches dans ce groupe. Ce critique soutient lui-même cette hypothèse dans ses publications [22, p. 569] et [23] :

> « *Derrière l'énumération à l'aide d'encoches se cache le nombre : l'histoire de l'arithmétique commence. On l'écrira un jour, grâce au grand nombre de monographies consacrées à la numération chez les primitifs, et on mettra en valeur les deux découvertes capitales que nous devons à nos ancêtres illettrés : le nombre et les systèmes de nombres, c'est-à-dire le nombre et le calcul.* »

8. « *… il est assez facile, moyennant quelque petits arrangements, de « faire parler » notre os, et même, en creusant un tout petit peu comme ci-dessus, de lui faire dire des choses contradictoires.* »

Bien sûr, et la théorie de Ramsey [24] en dit même plus sur toute série finie de nombres, mais restons-en aux « *petits arrangements* ». On pourrait y voir des nombres premiers (comme suggéré initialement par J. de Heinzelin), ou des triplets de Pythagore, etc. Cependant, si l'on prend la peine de considérer tous les groupements à la fois et de regarder non pas par le petit bout de la lorgnette mais de considérer l'ensemble, on ne peut qu'être amené à considérer une hypothèse arithmétique. J. de Heinzelin fut très prudent à ce sujet [21, pp 67-70] :

> « *… De l'avis des mathématiciens que j'ai consultés, aucun moyen logique ne peut prouver que ces chiffres sont dus ou non au genre de « hasard » qui intervient par exemple dans un compte de chasses ou de recettes. Cependant chacun oserait dire, je pense, dans l'explication de ce tableau de chiffres, que le sentiment des choses humaines fait pencher pour l'hypothèse arithmétique.*

> *S'il y a arithmétique, les calculs se fondent assurément sur les bases 2 et 10 ; l'usage primitif de celles-ci n'est pas fait pour étonner, car elles sont les plus naturelles à l'homme. … »*

Parmi les (neuf) scientifiques de renom de l'époque que de Heinzelin avait consultés figurent L. Hogben, l'auteur du livre *'Mathematics for the Million'* [25], et des spécialistes de l'Université Libre de Bruxelles et d'autres universités.

L'interprétation des nombres 11, 13, 17 et 19 d'encoches figurant sur la colonne de gauche du bâton comme un signe d'une connaissance des nombres premiers a été rejetée par les auteurs dans plusieurs publications. Cependant, la disposition des nombres 11, 13, 17 et 19 appelait une autre explication, par exemple en les considérant respectivement comme $12 \pm 1$ et $18 \pm 1$. On obtient alors des indications pointant vers la base 12 (et/ou ses sous-multiples 3, 4 et 6). De plus, les sommes des trois colonnes sont des multiples de 12. Les bases 10 et 20 sont physiologiquement évidentes, mais la base 12 l'est tout aussi bien. On compte avec le pouce d'une main les phalanges des quatre autres doigts, comme il est encore pratiqué de nos jours par certaines populations (voir [4]), et le nombre des douzaines sur l'autre main, donnant un total de 60. Cette interprétation suggère pourquoi 12 et 60 vont souvent de pair, même de nos jours, comme pour la subdivision du temps.

**Une règle à calcul ?**

9. *« Ils vont même plus loin que leur prédécesseur en affirmant que les nombres des trois colonnes sont liés de telle sorte que leur ensemble forme une règle à calcul. »*

Il n'a jamais été affirmé que l'os d'Ishango fût une règle à calcul. Les auteurs décrivent dans leurs publications l'arrangement et les particularités géométriques des encoches dans les différents groupements et émettent l'hypothèse que ce bâton aurait pu servir de support à une méthode de comptage, similaire à une règle à calcul, sans pour autant prétendre qu'il s'agissait d'un tel instrument. Par exemple [7, pp 342]

> *« … the proposed hypothesis of considering the bone as an ancient 'slide rule' to display simple addition arithmetic fits well with the various notch geometrical patterns. It shows a systematic occurrence of the derived base 12 and the central role that this base played in the proto-mathematics of the ancient Ishango people. »*

**Des additions inventées ?**

10. *« Il n'y a que quatre additions exactes ; mais comme les auteurs veulent que notre os soit une table d'additions, il leur faut en fabriquer d'autres de force. Par exemple, nous disent-ils, le 3 et le 6 du milieu sont presque en face du 11 de droite ; c'est donc que 3 et 6 ont été additionnés et le résultat mis à droite. Il est vrai qu'il manque 2 (entre parenthèses dans le tableau ci-dessus) : c'est donc que le 2 a été omis pour une raison inconnue ! C'est par le même procédé que Pletser et Huylebrouck inventent trois autres additions, notées sur la deuxième et les deux dernières lignes du tableau ci-dessus avec les nombres manquants entre parenthèses.*

*En supposant même qu'il s'agisse d'additions, quel profit tirer d'une table aussi brouillonne, où il faut prendre les nombres tantôt par deux, tantôt par trois, et où les résultats sont tantôt à droite et tantôt à gauche ?»*

Et plus loin : « *Au congrès de 2007, les auteurs fournissent un schéma supplémentaire en supposant que le cinquième nombre de la colonne du milieu est 10. Cette fois-ci, toutes les opérations sont fausses.* »

Les auteurs ont essayé d'envisager toutes les possibilités de combinaisons d'opérations d'arithmétique élémentaire entre les nombres d'encoches de la colonne du milieu et ceux des deux autres colonnes, tout en respectant l'arrangement et les formes des encoches [15]. Certaines tentatives donnaient un résultat, d'autres pas. Deux de ces tentatives ont été présentées dans ces tableaux, sans aucune intention de manipulation.

11. « *Et quel serait l'intérêt de telles additions, dans la mesure où la fusion de trois paquets de 3, 6 et 4 encoches, par exemple, en un seul paquet de 13 n'apporterait rien de plus qu'une lecture plus difficile de ce nombre ?* »

Si ce n'est qu'elle apporte une illustration visuelle de la somme des trois paquets de 3, 6 et 4 encoches de la colonne du milieu comme un nouveau paquet de 13 encoches dans la colonne de gauche.

**L'intention des hommes préhistoriques et inconnues**

12. « *Il est bien connu qu'avec les premières véritables notations numériques, une telle « addition » n'aurait eu aucun sens, puisque le nombre 13, pour reprendre le cas qui nous occupe, n'aurait jamais été écrit avec 13 encoches régulièrement espacées, mais uniquement sous la forme de paquets distincts, pour en faciliter la lecture ; en Égypte antique par exemple, le hiéroglyphe de 9 n'était pas 9 barres alignées à égale distance les unes des autres, mais 4 barres placées au-dessous de 5 autres, ou plus souvent trois paquets de 3 barres les uns en dessous des autres.* »

Est-ce une critique ou un argument ? Les auteurs ne comprennent plus ce critique, car d'une part, il décrie leurs interprétations similaires, et d'autre part, voici qu'il propose lui-même l'idée audacieuse du passage d'une représentation d'un paquet de 13 encoches à une représentation plus abstraite d'un nombre, 13. Et là, n'était certainement pas l'intention des hommes préhistoriques d'Ishango qui n'avaient sans doute pas encore découvert la notion des nombres et de leur représentation. Ce passage à l'abstraction et à la représentation des nombres est bien postérieur à la civilisation préhistorique d'Ishango.

De plus, ce critique signale dans sa note 9 de bas de page 9 de [19], la numérotation babylonienne utilisant un symbole (le clou) en base 10 et un autre (le chevron) en base 6. Le nombre 9 s'écrit ainsi avec neuf clous, comme le signale B. Rittaud [26, p.2], cité par ce critique. Alors, deux poids, deux mesures ?

13. « *L'affaire ne marche pas, donc il y a des intentions inconnues ! À cela s'ajoutent des spéculations peu convaincantes sur des comparaisons des longueurs ou des inclinaisons des encoches, qui de toute manière ne parviennent pas à justifier les ajouts nécessaires de 2, 1 ou -1 pour donner un semblant de cohérence à l'ensemble.* »

Concernant les « *intentions inconnues* » de la première phrase, les auteurs signalent qu'il arrive à ce critique de confesser lui-même son ignorance, par exemple [27, p. 73]

> « *Pour des raisons que l'auteur ignore, les symboles primitifs sont essentiellement des symboles géométriques* »

De plus, ces « *spéculations peu convaincantes* » font allusion à l'étude de la forme et de la taille des encoches des différentes colonnes qui sont détaillées dans [15].

**Le choix des nombres sur l'os**

14. « *S'agissant du cinquième nombre de la colonne du milieu, on peut se demander pourquoi les auteurs ont choisi de prendre 10, avec lequel toutes les additions sont fausses, plutôt que 9, avec lequel quatre sont exactes. La raison est qu'avec 10, le total des nombres de la colonne du milieu est 48, qui est un multiple de 12, comme le total de 60 des colonnes de droite et de gauche.* »

Effectivement.

15. « *Pour rendre compte du choix des nombres présents sur l'os, Pletser et Huylebrouck avancent en effet l'idée que :*
> « *Les nombres 3 et 4 pourraient avoir constitué la base du système arithmétique en usage dans l'ancienne population d'Ishango pour opérer sur les petits nombres, et la base dérivée 12 pour les grands nombres.* »

*Où voit-on les bases 3 et 4 ? Dans la colonne du milieu, selon les auteurs, puisque nous avons, de haut en bas :*
- *3 puis 6, donc 3 puis 3×2*
- *4 puis 8, donc 4 puis 4×2*
- *9 ou 10, donc 4×2+1 ou 4×2+2*
- *deux fois 5 « pour montrer deux manières d'obtenir le nombre composé 5 en ajoutant 1 ou 2 aux bases 3 et 4 »*
- *7 « montre comment obtenir le nombre composé 7 en additionnant les bases 3 et 4».* »

Effectivement.

**La combinaison de plusieurs bases**

16. « *Mais si 5 est écrit deux fois parce qu'il y a deux bases, pourquoi les autres nombres 6, 8, 9 ou 10 et 7 ne figurent-ils qu'une seule fois ? Si l'objectif, d'autre part, avait été de mettre une base en relief, celle-ci apparaîtrait nettement ; on devrait voir clairement deux ensembles de 3 au sein du paquet de 6, deux ensembles de 4 au sein du paquet de 8 et ainsi de suite. Or il n'est rien de tel, il n'y a aucun groupement visible régulier faisant penser à une base.* »

Il n'a été dit nulle part que le sculpteur d'encoches avait « *l'objectif de mettre une base en relief* » sur cet os. Les auteurs ne comprennent pas d'où ce critique sort cette remarque.

Les auteurs ont avancé que ce bâton est probablement le témoignage d'un peuple qui comptait en mélangeant les bases 10 et 12 (ou 6). Un européen écrirait peut-être ||||| ||||| |||| pour compter 14 jours, un homme d'Ishango |||||   |||||||   || . La combinaison de plusieurs bases est chose courante et n'est pas étonnante. Par exemple, en français, 77 et 93 réfèrent bien à un passé à bases 10 et 20. Si dans 22.000 ans on retrouve un texte français avec les expressions 'soixante-dix-sept' et 'quatre-vingt-treize', l'archéologue du futur pourrait en conclure avec raison que la langue française mélangeait les bases 10 et 20. De plus, il y aurait peu de chances que ce texte soit une dissertation savante sur l'arithmétique expliquant l'utilisation de ces différentes bases.

Ensuite, la duplication n'est pas nécessairement une opération d'ajout à l'original d'une copie de l'original, ni une simple multiplication par 2. Il peut s'agir d'une reconstruction d'un nouveau groupe ayant un nombre d'encoches double de celui du groupe initial mais arrangé différemment [15]. Le fait d'exiger de « *voir clairement deux ensembles de 3 au sein du paquet de 6, deux ensembles de 4 au sein du paquet de 8 et ainsi de suite* » est ce qui pourrait être appelé 'regarder par le petit bout de la lorgnette' et fait preuve d'un égocentrisme mathématico-culturel moderne trop simpliste.

17. « *Où voit-on la base 12 ? D'une part, comme nous l'avons dit, dans les totaux de chaque colonne qui sont des multiples de 12. Et d'autre part, affirment les auteurs, par le fait que dans la colonne du milieu :*
   - *6 intervient dans deux des additions supposées (deux premières lignes du tableau ci-dessus) : 2 ×6 = 12*
   - *4 intervient dans trois des additions supposées : 3 ×4 = 12*
   - *8 intervient dans trois des additions supposées : 3 ×8 = 24 = 2 ×12* »

Effectivement, mais, comme dit plus haut, la base 12 apparaît également dans la colonne de gauche avec les nombres 11, 13, 17, 19, qui après l'abandon de l'hypothèse des nombres premiers, sont vus comme $12 \pm 1$ et $18 \pm 1$.

18. « *Nous avons donc la situation suivante : 6 n'est gravé qu'une fois, sous forme d'un groupe de six encoches de la colonne du milieu ; mais comme on a supposé qu'il intervient dans deux additions, et bien que celles-ci soient fausses à cause d'une intention inconnue, cela donne 12 ! De même pour 4 et 8 qui interviendraient trois fois chacun, donnant respectivement 12 et 2 ×12. Est-il possible d'être convaincu par de tels tours de passe-passe dignes de la plus plate littérature numérologique ?* »

Remettons les choses à leur place. Les bases 3 et 4 sont déduites des opérations élémentaires dans la colonne du milieu et la base 12 est déduite des nombres $12 \pm 1$ et $18 \pm 1$ de la colonne de gauche. Les auteurs constatent simplement que les sommes des nombres de la colonne du milieu sont égales à 12 ou à des multiples ou sous-multiples de 12. Ce fait n'est nulle part utilisé comme argument pour renforcer la véracité de l'hypothèse de la base 12 et de ses sous-multiples.

**Le deuxième os**

19. « *Jean de Heinzelin a trouvé en 1959, toujours à Ishango, un deuxième os marqué d'encoches, et a proposé en 1998 une interprétation, qui, d'après Pletser et Huylebrouck, confirme ce qui précède. Il ne sera pas utile d'ennuyer beaucoup plus longtemps le lecteur avec cela.* »

Les auteurs doutent que ce critique ait vu le second bâton, car rien n'avait été publié à ce sujet avant 2007, à l'occasion du congrès auquel ce critique était convié, mais qu'il a préféré éviter.

Les auteurs pensent que ce point au contraire est du plus grand intérêt. Sur son lit de mort, Jean de Heinzelin avait changé d'avis au sujet de l'interprétation des encoches du premier bâton, y voyant effectivement un outil arithmétique élémentaire utilisant plusieurs bases simples. Cette conclusion n'avait pas été publiée suivant ses dernières volontés et ne pourrait être publiée que quelque temps après sa mort. La surprise des auteurs fut grande de voir leur hypothèse confirmée a posteriori par de Heinzelin, en même temps que l'annonce de l'existence du second os.

20. *« Il sera suffisamment édifié en jetant simplement un coup d'œil sur la figure 6 et en lisant ce qui suit :*

> *« De Heinzelin ajoutait que le petit tiret sur la colonne E était à la dixième place, et se demandait si cela annonçait « un passage de la base 10 à la base 12 » […] Puisque la colonne C a un total de 20 encoches, et la colonne E 18, les bases 6 et 10-20 semblent se faire jour. De plus, il y a deux concordances spatiales entre les rangées : E10 = F1 = G10 et E12 = F2 = G12. »*

*Cet os, selon des hypothèses de De Heinzelin reprises par Pletser et Huylebrouck, pourrait être le témoignage d'un changement de base, avoir eu une fonction didactique, ou même avoir servi dans l'échange entre groupes ethniques, l'un pratiquant la base 10, l'autre une autre base comme 12, 16 ou d'autres encore. »*

L'archéologue de Heinzelin fut une sommité mondiale de l'archéologie. Ses arguments étaient toujours fondés et avaient du poids.

**Le petit bout de la lorgnette**

21. *« Rappelons encore une fois que l'on ne peut sérieusement soutenir la présence de bases que si l'on est en face de regroupements clairs et systématiques ; il ne suffit pas que 18 = 3×6 pour démontrer pas la présence d'une base 6 ! »*

Et rappelons également que cette manière de voir les choses est contraire à l'esprit scientifique qui présuppose une ouverture d'esprit et l'examen de toutes les hypothèses possibles et plausibles. Ne chercher que « *des regroupements clairs et systématiques* » est tellement réducteur et fait preuve d'une étroitesse de vue (le petit bout de la lorgnette…) passant toute observation au moule de notre vision moderne et actuelle des choses sans s'interroger sur la manière dont un homme préhistorique aurait pu voir les choses de son point de vue.

Effectivement, 18 = 3 x 6 n'est pas suffisant pour en déduire l'existence d'une base 6 s'il n'y a qu'une seule indication de ce type. Par contre, si plusieurs indications sont présentes comme montré par de Heinzelin et d'autres scientifiques pour le premier bâton d'Ishango avec plusieurs exemples d'opérations élémentaires de duplication et d'addition faisant intervenir les nombres 3, 4 et 12, à côté du nombre 10, tout chercheur ayant été formé à la méthode scientifique, c.-à-d. l'observation et la déduction, est en droit de se poser la question s'il n'existe pas effectivement plus qu'un simple hasard et si des corrélations entre les différentes bases ne devraient pas être considérées.

**Le comparatisme ethnographique**

22. *« Les quelques exemples ethnographiques donnés par les auteurs dans les Actes du Congrès de 2007 ne peuvent pas davantage nous convaincre. Car si, chez tel peuple du Congo, on dit l'équivalent de 'douze-un' pour notre treize, on peut bien parler de base 12, mais cela signifie justement que la seule façon de dire 13 est de prononcer 12, puis 1 ; l'équivalent graphique serait de graver un paquet de 12, puis une autre encoche bien séparée. Il en est de même avec les nombres simplement figurés avec les doigts. Les Shambaa de Tanzanie montrent 6 en étendant trois doigts de chaque main, et disent l'équivalent de 'trois-trois' pour six. Même accord entre le geste et le mot pour 8 qui se dit 'quatre-quatre' et se montre avec quatre doigts de chaque main. Il y a désaccord pour 7, qui se dit 'dix moins trois'* *et se montre par 4 doigts de la main droite et 3 de la main gauche ; mais le geste de 7, comme celui de 8 et de 6, sépare nettement les bases 3 et 4, si bases il y a. Or de telles séparations claires systématiques en sous-ensembles de 3, 4 ou 12 encoches n'existent sur aucun des os d'Ishango. Les exemples ethnographiques ne font donc qu'enfoncer un peu plus les théories de Pletser et Huylebrouck. »*

Il est surprenant de lire le début de cet extrait sous la plume de quelqu'un qui se dit le chantre du comparatisme ethnographique, en s'extasiant même [22, pp. 563-564] :

> *« Le comparatisme ethnographique, qui reconnait cette analogie, ouvre une voie de recherche à peine explorée pour les mathématiques de la préhistoire, mais riche de promesses puisqu'il permet, dans une certaine mesure, de faire parler les trouvailles archéologiques. Et il impose aussi aux hypothèses et théories fondées uniquement sur les documents préhistoriques d'être vérifiées par les documents ethnographiques et inversement. »*

Les auteurs ne voient pas comment ni pourquoi « *Les exemples ethnographiques ne font donc qu'enfoncer un peu plus les théories* » des auteurs. Rien ne permet d'affirmer cela comme le fait ce critique. Au contraire même, si celui-ci s'était donné la peine de lire avec attention les différents arguments circonstanciels ethnographiques donnés dans [15].

**Une question de femmes battues ?**

23. *« Néanmoins, que les prétendues lectures soient simples ou sophistiquées, toutes reposent sur un même fondement arbitraire, qui consiste à croire que des encoches ont nécessairement un caractère numérique. Claudia Zaslavsky raconte que certaines femmes africaines font de temps en temps une encoche dans le manche de leur cuillère en bois. Marquent-elles des jours ? Jouent-elles avec des nombres ? Nullement : elles font une marque chaque fois qu'elles reçoivent un coup de leur mari ; et dès que le manche de la cuillère est rempli, elles demandent le divorce. »*

Les auteurs ont cité Zaslavsky dans leur référence 23, p. 166 de [15]. Plus même, ils ont cité des rapports semblables écrits par des missionnaires et administrateurs du Congo belge qui décrivaient comment et pourquoi de tels bâtons à encoches étaient utilisés.

De plus, il est exact que cette hypothèse ne peut pas non plus être entièrement rejetée. Mais pourquoi alors une femme préhistorique d'Ishango aurait fait des marques aussi sophistiquées avec des groupements qui présentent une telle régularité d'ensemble ? Et surtout pourquoi faire des marques sur trois colonnes ? On peut imaginer qu'elle aurait pu déjà demander le divorce après avoir rempli la colonne du milieu avec 48 encoches, ou à la rigueur celle de gauche ou celle de droite, avec 60 encoches chaque fois !

Bien que ne pouvant pas être formellement rejetée, il est très probable que cette hypothèse serait la moins plausible de toutes pour le premier bâton d'Ishango.

24. *« Une encoche peut donc n'être qu'une marque, ce qui ne paraît pas grand-chose si l'on est obnubilé par l'arithmétique ; il y a pourtant là la plus importante invention que nous devons à nos ancêtres du Paléolithique supérieur, celle du signe. »*

Effectivement.

25. « *Et à s'égarer dans des spéculations numériques hasardeuses, on gaspille du temps, de l'argent et du papier alors qu'il y a tant à découvrir dans les signes préhistoriques, y compris en ce qui concerne la gestation intellectuelle du concept de nombre, en les confrontant avec la documentation ethnographique.* »

Si ce critique avait suivi cette recommandation à la lettre, son œuvre scientifique, se limitant à quelques publications au niveau national et répétant les mêmes critiques que celles discutées ici, aurait disparu entièrement.

26. Les deux contre-exemples ethnographiques donnés par ce critique [19, pp.14-17] sont excellents mais ne sont pas des contre-exemples, bien au contraire. Les auteurs ont repris des exemples semblables dans certaines de leurs publications. Ces histoires montrent surtout un point important : la connaissance a priori de la raison d'être de ces supports marqués ou signés. Autrement dit, une autre source d'information (un écrit, une tradition orale, une information supplémentaire, quel que soit le support) permet le déchiffrement aisé de ces objets.

Les os d'Ishango ont été retrouvés avec d'autres outils et des pointes de harpon [21]. Aucun de ces autres artefacts n'indique un quelconque lien avec les encoches sur les os. Le critique le dit d'ailleurs lui-même [22, p. 567] en parlant de l'os d'Ishango :

> « *De la même façon, aucun document ethnographique ne permet d'étayer la thèse mentionnée plus haut […] d'un chasseur-cueilleur fabriquant une table de nombres premiers, une table de doubles, …* »

Alors que faire ? Ignorer ces artefacts et les laisser dormir dans un tiroir ou au contraire les étudier et essayer de les déchiffrer ?

Rappelons que l'étude du premier os a été réalisée sans même savoir qu'un second os avait été découvert au même endroit et que la description de ce second os n'a été publiée que plusieurs années après la mort de son découvreur. Or les conclusions auxquelles sont arrivées les auteurs sur le premier os ont été corroborées par après par les conclusions de Jean de Heinzelin sur son lit de mort avant la publication des résultats des auteurs et sans que ceux-ci n'en aient connaissance.

**Quelle calculette ?**

27. « *S'il est désolant que des fictions du style « calculette d'Ishango » circulent et soient souvent prises pour argent comptant, ce n'est pas seulement en raison de leur faiblesse interne et de leur invraisemblance, mais c'est aussi qu'il y aurait mieux à faire avec toute cette documentation venue de l'archéologie et de l'ethnographie.* »

Les auteurs n'ont jamais employé l'expression « *calculette d'Ishango* » dans leurs publications et communications. Ils laissent à ce critique la paternité de cette expression.

## **Conclusions:**

**Une conclusion correcte**

L'étude du premier os est bien arrivée à une conclusion correcte, bien que peu spectaculaire. Elle a été réalisée sans savoir qu'un second os avait été découvert au même endroit et avant la publication de sa description plusieurs années après la mort de son découvreur. Or les conclusions auxquelles sont arrivées les auteurs sur le premier os ont été corroborées par après par les conclusions de Jean de Heinzelin sur son lit de mort avant la publication des résultats des auteurs et sans que ceux-ci n'en aient connaissance.

**Un défaut épistémologique**

Une meilleure recherche bibliographique aurait permis à ce critique de se faire une opinion sur l'ensemble des éléments, mais on peut douter que ce critique ait une maîtrise suffisante de la langue anglaise, car il tire des éléments hors de leur contexte et utilise des sources de seconde main au lieu de consulter les documents originaux. Enfin, il est vrai que ce n'est pas la première fois que ce critique a des problèmes à lire des publications scientifiques et à se documenter proprement (Voir par exemple [28]).

**Non-respect de droits d'auteurs**

Les auteurs signalent également l'utilisation par ce critique d'images sans respect du droit d'auteur ; ses illustrations sont reprises d'autres articles sans aucune demande d'autorisation, ou du site du Musée des Sciences Naturelles de Bruxelles, sans aucun respect des lois en vigueur.

**Un vocabulaire douteux, peu académique**

Finalement, les auteurs ne résistent pas au plaisir de relever dans ce petit pamphlet les adjectifs et expressions suivantes qui éclairent sur le style de ce critique : « *ridicule»* ([19], p. 5), « *en fabriquer d'autres de force»* (p. 7), « *aussi brouillonne»* (p. 7), « *L'affaire ne marche pas»* (p. 8), « *des spéculations peu convaincantes»* (p. 8), « *un semblant de cohérence»* (p. 8), «*de tels tours de passe-passe dignes de la plus plate littérature numérologique»* (p. 10), « *le même « tour »»* (p. 10), « *il faut fabriquer des 12»* (p. 10), « *des bricolages aussi confus»* (p. 10), « *Il ne sera pas utile d'ennuyer beaucoup plus longtemps le lecteur avec cela. Il sera suffisamment édifié »* (p. 10), « *On invente ensuite quelque contexte concret et on imagine une histoire pour ne fabriquer en fin de compte qu'une fiction mathématique. C'est ce à quoi nous venons d'assister avec les spéculations de Heinzelin, Pletser et Huylebrouck, lesquels ne sont d'ailleurs que les derniers en date d'une longue lignée de victimes de l'illusion mathématicienne ; illusion d'autant plus tentante et fréquente qu'il s'agit de préhistoire. »* (p. 12), « *Les sirènes de l'illusion mathématicienne »* (p. 12), « *le mathématicien naïf »* (pp. 12-13), « *Quatre est féminin parce que les femmes ont quatre lèvres. »* (p. 15), « *suffisamment cruels pour les auteurs de fables mathématiques»* (p. 17).

Quant aux termes « *faiblesse interne et […] invraisemblance »* (p. 17), ils seraient mieux accolés à ce petit pamphlet critique.

Les auteurs déplorent ce langage peu académique mais ils ne veulent pas entrer dans une vaine polémique pour les raisons expliquées dans l'introduction.

Cependant, ils sont ouverts aux débats d'idées et ils seraient heureux de participer à un échange de points de vue lors d'une conférence.

*V. Pletser, D. Huylebrouck, août 2015*

Références

[1] D. Huylebrouck, 'Puzzles, patterns, drums: the dawn of mathematics in Rwanda and Burundi', *Humanistic Mathematics Network*, Journal no 14, November 1996.

[2] D. Huylebrouck, 'The Bone that began the Space Odyssey', *The Mathematical Intelligencer*, Vol. 18, no 4, p. 56, Autumn 1996.

[3] D. Huylebrouck, 'Histoires célestes du Rwanda', *Ciel et Terre, Bulletin de la Soc. Royale Belge d'Astronomie, Météorologie et Physique du Globe*, Vol. 113 (2), 83-89, 1997.

[4] D. Huylebrouck, 'Counting on hands in Africa and the origin of the duodecimal system', *Wiskunde en Onderwijs (Mathématiques et Education),* no 89, 1997.

[5] V. Pletser, D. Huylebrouck, 'The Ishango Artefact: the Missing Base 12 Link.' *Proc. Katachi Univ. Symmetry Congress (KUS2),* T. Ogawa, S. Mitamura, D. Nagy & R. Takaki (ed.), Paper C11, Tsukuba Univ., Japan, 18 Nov. 1999.

[6] V. Pletser, D. Huylebrouck, 'Research and promotion about the first mathematical artefact: the Ishango bone', *Proc. PACOM 2000 Meeting Ethnomathematics and History of Mathematics in Africa*, Cape Town, South Africa, 1999.

[7] V. Pletser, D. Huylebrouck, 'The Ishango Artefact: the Missing Base 12 Link', *Forma*, 14-4, 339-346, 1999.

[8] D. Huylebrouck, 'Afrika + Wiskunde' (in Dutch), VUBPress of the Free University of Brussels, 2005; traduit en français 'Afrique + Mathématiques' VUBPress of the Free University of Brussels, 2008.

[9] D. Huylebrouck, 'De parabolische vlucht van het Ishangoartefact, de oudste wiskundige vondst' (Le vol parabolique de l'artefact d'Ishango, la plus ancienne découverte mathématique), *Wiskunde en Onderwijs (Mathématiques et Education),* no 101, 2000.

[10] D. Huylebrouck, F. Dumortier, 'A Stamp for World Mathematical Year', *European Mathematical Society Newsletter*, p. 27, November 2000.

[11] D. Huylebrouck, 'Middle column of marks found on the oldest object with logical carvings, the 22000-year-old Ishango bone from the Congo', sequence A100000, *On-Line Encyclopaedia of Integer Sequences*, http://oeis.org/, November 2004

[12] D. Huylebrouck, 'L'Afrique, berceau des mathématiques', *Pour la Science*, Dossier 47, avril-juin 2005 (6 pages). L'article fut traduit en portugais (au Brésil) et en allemand.

[13] D. Huylebrouck, 'Mathematics in (central) Africa before colonisation', *Anthropologica et Praehistorica, Bulletin of the Royal Belgian Association for Anthropology and Prehistory*, 135 – 162, Vol. 117, 2006.

[14] D. Huylebrouck, 'Report: The ISShango project', *Journal for Mathematics and the Arts*, Taylor and Francis, September 2008.

[15] V. Pletser, D. Huylebrouck, 'An interpretation of the Ishango rods', *Proc. Conf. 'Ishango, 22000 and 50 years later: the cradle of Mathematics?'*, D. Huylebrouck ed., Royal Flemish Academy of Belgium, KVAB, 139-170, 2008.

[16] D. Huylebrouck, 'Foreword to the book Lusona: find the omitted drawings' by Paulus Gerdes (Mozambique), Lulu Editions, July 2009.

[17] D. Huylebrouck, R. Matteus Berr, 'Vermessung', in '*Kunstlerhaus catalogue EVO EVO! 200 Jahre Charles Darwin'*, I. and P. Braunsteiner eds, Künstlerhaus ISBN 978-3-900926-84-7 Verlag Lehner Wien, 2009.

[18] D. Huylebrouck, 'België + wiskunde' (Belgique + Mathématique), Academia Press Gent, Belgique, 2013.

[19] O. Keller, 'Les fables d'Ishango, ou l'irrésistible tentation de la mathématique-fiction', *Bibnum*, https://www.bibnum.education.fr/sites/default/files/Ishango-analyse.pdf, août 2010.

[20] J. Yellen, 'Ishango and its importance for later African prehistory' (Ishango et son importance pour la préhistoire africaine), *Proc. Conf. 'Ishango, 22000 and 50 years later: the cradle of Mathematics?'*, D. Huylebrouck ed., Royal Flemish Academy of Belgium, KVAB, 63-81, 2008.

[21] J. de Heinzelin de Braucourt, 'Exploration du Parc National Albert', Fascicule 2, Les Fouilles d'Ishango, Instituts des Parc Nationaux du Congo Belge, Bruxelles, 1957.

[22] O. Keller, 'Questions Ethnographiques et Mathématiques de la Préhistoire', *Revue de synthèse* : 4$^e$ S. n° 4, oct.-dec. 1998, p. 545-573. Disponible sur http://rd.springer.com/article/10.1007/BF03181393

[23] O. Keller, « Préhistoire de l'arithmétique, La découverte du nombre et du calcul », in 'Si le nombre m'était conté', I.R.E.M. - Histoire des mathématiques, Paris, Ellipse, 2000.

[24] F.P. Ramsey, 'On a Problem of Formal Logic', *Proceedings London Mathematical Society*, s2-30 (1): 1930, 264–286. doi10.1112/plms/s2-30.1.264

[25] L. Hogben, « Mathematics for the Million», London, Allen & Unwin, 1936. ISBN-13: 978-0393310719

[26] B. Rittaud, 'À un mathématicien inconnu !', *Bibnum*, https://www.bibnum.education.fr/sites/default/files/RITTAUD_YBC7289.pdf, dernier accès 09/08/2015.

[27] O. Keller, « Préhistoire de la Géométrie : Premiers éléments d'enquête, premières conclusions », *Science et Techniques en Perspective*, Vol. 33, Université de Nantes, 1995.

[28] J. Høyrup, 'Book review of 'Préhistoire de la Géométrie : Premiers éléments d'enquête, premières conclusions', published in *Isis* 87,1 996, 713-714. Disponible sur http://www.academia.edu/3131829/Prehistoire_de_la_Geometrie_Premiers_elements_denquete_premieres_conclusions
http://akira.ruc.dk/~jensh/publications/1996%7BR%7D17_Keller_Prehistoire_MS.PDF

# Contradictions and narrowness of views in "The fables of Ishango, or the irresistible temptation of mathematical fiction", answers and updates

Vladimir Pletser ([1]), Dirk Huylebrouck ([2])

*(1) Sciences Department, Manned Spaceflight and Operations Directorate, European Space Research and Technology Centre, European Space Agency, Noordwijk, The Netherlands;*
*Present address: Technology and Engineering Centre for Space Utilization,*
*Chinese Academy of Sciences, Beijing, China*
*(2) Faculty of Architecture, University of Leuven, KUL Brussels-Gent, Belgium*

Note : This is a translation of the original text in French titled "*Contradictions et étroitesse de vues dans « Fables d'Ishango, ou l'irrésistible tentation de la mathématique-fiction », réponses et mises au point*", June 2016. Available at http://arxiv.org/abs/1607.00860

**Preliminary**: *Response to a rumour: what is it?*

*A discovery made in the middle of the 20th century in the East of the Democratic Republic of Congo, in the heart of the Virunga forest, a forest declared heritage of humanity by UNESCO, of a bone with notches dating from more than 20,000 years ago that most likely influenced our science and that became a symbol.*

*The authors of the following lines, like many specialists with undoubtable qualifications, have issued their hypotheses. The French "Le Robert" dictionary defines the word hypothesis in a mathematical sense: "proposal admitted as a problem data" and in common parlance: "conjecture about the explanation or the possibility of an event". It specifies that the opposite is equivalent to: "conclusion, certainty, obvious".*

*Each discovery has a birth certificate established by specialized publications giving rise to international scientific conferences where authorized experts give their opinion. Its impact is evolving and assessed in this way.*

*But now, alongside these argued opinions, a denigrating rumour circulates. Mr. Olivier Keller, at ease in spreading rumours and using a loose vocabulary, acquires a certain notoriety in discrediting in a blog the thought of scientists who are more nuanced and cautious by trade. Should one ignore it because, always absent in the peer reviewed publications and appropriate conferences, no immediate response and no direct debate with the author of this blog were possible?*

*We asked the present authors to answer these criticisms. They chose to track closely each of the 28 allegations of this blog, to answer them point by point.*

*Gabriel Castel, August 2015*
*Committee Ishango Milele*
*www.ishango-milele.com*

## Introduction

**The hypotheses**

Since the end of the 1990s, the authors have introduced and explained the following hypotheses about the Ishango bone, discovered by archaeologist Jean de Heinzelin in the late 1950s:

1) The Ishango bone, dated from 22,000 years, can be considered as the oldest mathematical tool of humankind because the arrangement of the notches on three columns suggests an arithmetical intention.
2) In addition, it appears that several bases are used in this elementary arithmetic: the base 10 and base 12 with its submultiples 3, 4 and 6. The geometric arrangement of the notches in the various groupings on the three columns allows to compute other basic arithmetic operations.

These hypotheses were presented in several publications [1-18] and in international conferences on mathematics and ethnomathematics.

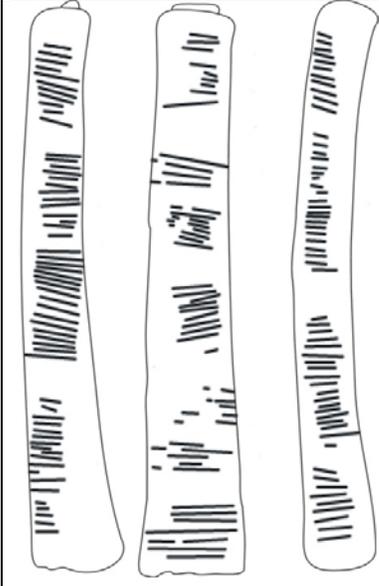

*Schematics of the first Ishango bone (adapted from Wikimedia Commons) and the corresponding numbers of notches of the three columns.*

The authors have explored in detail the arithmetic relationships between the numbers of notches. One knows for example that the sums of the numbers in the columns from left to right are 60, 48 and 60, that the numbers in the left column can be seen as primes, that doublings are seen in the middle column, and that the numbers in the right column can be interpreted as 10 and 20 plus or minus 1.

**The conclusion.**

The authors have proposed various hypotheses, and rejected others (such as the one of prime numbers), yielding an unspectacular conclusion: it is probably a tool that counts events or things, noted by someone who mixed bases 10 and 12.

However, the authors did not enter into discussions about the use of this rod: is it an arithmetic game, the hypothesis of the archaeologist Jean de Heinzelin; or a calendar, according to Marshack; or a simple counting tool; or something else? The authors have left this discussion to ethnologists and anthropologists.

**A loose academic criticism**

Olivier Keller, in his single publication [19], has criticized, even ridiculed, the work of the authors concerning this rod. It took them some time before they decide to respond to these criticisms. Indeed, the work of this critic has never been published in international journals or journals with peer reviewing. In addition, this critic has never participated in any international scientific congress, where he could have defended his point of view, which is part of the scientific process.

Nowadays, it is impossible to respond to all the criticisms on blogs or websites flourishing among certain pressure groups or defending some fanatical or political beliefs. Without wanting to polemicize, the authors want to present the following few arguments to refute these criticisms one by one, in italics below and reproduced as such from [19].

## **The answers**

**The bone, the tool and its age**

1. *« … the first of the two Ishango bones – rose to fame by being presented as a noteworthy scientific text … »*

Strange formulation: a bone is not a *«scientific text »*.

2. *« A fragment of quartz affixed to one end shows that it was a tool handle; … »*

This hypothesis is supported by renowned archaeologists who studied the bone: Jean de Heinzelin (B), Alison Brooks (USA), John Yellen (USA), Els Cornelissen (B), and is also mentioned on the website of the Institute of Natural Sciences in Brussels where the object is exposed.

3. *« … It is usually dated to 20,000 years BCE»*

The accepted dating is of 22,000 years. This dating is obtained by the carbon-14 method, and confirmed by several other archaeological methods, as explained by Alison Brooks. The Ishango civilization might even be 90 000 years old, following the work of J. Yellen [20].

**The suspicion on the notches**

4. *« Several of the notches are worn way or barely visible, which immediately makes any numeric interpretation suspect. »*

This critic should be advised to not rely on the observation of photographs or reproductions, but to consult the original work of J. de Heinzelin [21, pp 64-70] describing the circumstances of its discovery, unless of course one should also suspect the illustrious archaeologist of intellectual dishonesty. Moreover, the technician and collaborator of de Heinzelin, Marcel Spinglaer, was a recognized specialist who can certainly not be suspected.

5. *« Let's take the middle column: according to the author, 3 is doubled to 6, 4 to 8 and 5 to 10. But the 5 and the 10 are doubtful: one of the sets of 5 is genuinely illegible, and in reality the 10 could be a 9. In addition, in the case of a duplication of 5, 3 and 4, there is to explanation* *as to why the set of five notches is shown twice, whereas the group of three and five* [sic] *are shown only once. »*

Where is the problem? Couldn't we see just two operations of duplication (3 => 6 and 4 => 8) and an addition operation (10 = 5 + 5)?

6. *« And what is the role of 7, which is neither involved in duplication nor doubled? Unless the bottom of the middle column reads 10, 4, 5 and 7 (and not 10, 5, 5 and 7), which would give us 7 doubled, with 10+4, and 5 doubled, with 10. »*

Isn't the author of these lines engaging in speculation that he himself so strongly decries? Nevertheless, the proposed interpretation is interesting and we leave the responsibility to its author.

7. *« Once it has been decided that the sets of notches are numbers, … »*

It is true that the starting hypothesis of any mathematical interpretation of this Ishango rod is the association with a group of notches of the number of notches in this group. This critic himself supports this hypothesis in his publications [22, p. 569] and [23]:

> *« Behind the enumeration using notches hides the number: the history of arithmetic begins. One will write it someday, thanks to the large number of monographs devoted to the numeration of the primitives, and one will highlight two main discoveries that we owe to our illiterate ancestors: the number and the systems of numbers, i.e. the number and the calculation. »*

8. *« … it's easy – given a few arrangements here and there – to load the bone with meaning, or even, if one pursues the argument a little further, as above, to make it say contradictory things. »*

Of course, and Ramsey theory [24] tells even more about any series of numbers, but let's stick to the *«few arrangements here and there»*. One could see here primes (as originally suggested by de Heinzelin), or Pythagorean triplets, etc. However, if we care to consider all groupings at once and look, not by the small end of the spyglass[1], but at the whole, one can only be led to consider an arithmetic hypothesis. J. de Heinzelin was very careful about this [21, pp 67-70]:

> *«… In the view of the mathematicians I consulted, no logical way can prove that these numbers are due or not to the kind of 'chance' that occurs for example in a hunting or revenue account. However, everyone would dare to say, I think, in the explanation of this table of numbers, that the feeling of human things is tipped for the arithmetic hypothesis. If there is arithmetic, calculations are certainly based on bases 2 and 10; the primitive use of these is not surprising, because they are the most natural to man. »*

The nine renown scientists of the time whom de Heinzelin had consulted include L. Hogben, author of the book *Mathematics for the Million* [25], and specialists of the Free University of Brussels, and of other universities.

The interpretation of the numbers 11, 13, 17 and 19 of notches on the left column of the rod as a sign of a knowledge of prime numbers was rejected by the authors in several publications.

---

[1] The French expression stating "*regarder par le petit bout de la lorgnette*" (translating literally as 'looking by the small end of the spyglass') means in fact not seeing the whole, but concentrating on a tiny part and missing the big picture.

However, the disposition of the numbers 11, 13, 17 and 19 called for another explanation, for example by considering them respectively as 12 ± 1 and 18 ± 1. One then get indications pointing to the base 12 (and/or its submultiples 3, 4 and 6). In addition, the sums of three columns are multiples of 12. Bases 10 and 20 are obvious considering human anatomy of hands and feet, but the base 12 is just as well. One counts with the thumb of one hand the phalanges of the four fingers, as it is still practiced today by some populations, and the number of dozens on the other hand, giving a total of 60. This interpretation suggests why 12 and 60 often go together, even today, as for the time subdivision

**A slide rule?**

9. *« Indeed, they go further than their predecessor by affirming that the numbers in the three columns are related in such a way that together they form a calculation rule* [sic]. *»*

It was never stated that the Ishango bone was a slide rule. The authors describe in their publications the arrangement and the geometrical characteristics of the notches in the various groupings and hypothesize that this rod could have served as support for a method of counting, similar to a slide rule, without being necessarily such an instrument. For example, [7, pp 342]:

> *« […] the proposed hypothesis of considering the bone as an ancient 'slide rule' to display simple addition arithmetic fits well with the various notch geometrical patterns. It shows a systematic occurrence of the derived base 12 and the central role that this base played in the proto-mathematics of the ancient Ishango people. »*

**Invented additions?**

10. *« Only four of the additions are exact. However, as the authors want this bone to be an addition chart, they have to forcibly make up others. For example, the 3 and the 6 in the middle column, they tell us, are almost aligned with the 11 in the right column, ergo the 3 and the 6 have been added together and the answer shown on the right. True, the answer is out by two (cf. the +2 in the table above), ergo the 2 has been left out for some unknown reason! Pletser and Huylebrouck use the same technique to invent three other additions, shown in the second and last two lines of the table above, with the missing numbers in parentheses. Let's suppose for a moment that this is an addition table. What is the point of such a muddled table, whose numbers have to be considered now in sets of two, now in sets of three, and whose answers are sometimes shown on the right and sometimes on the left? »*

And further: *« At the conference in 2007, the authors provided an additional chart and assumed that the fifth number in the middle column is 10. This time, all the operations are incorrect. »*

The authors tried to consider all possibilities of combinations of operations of elementary arithmetic between the notch numbers of the middle column and those of the other two columns, while respecting the arrangement and shapes of the notches [15]. Some attempts gave a result, others did not. Two of these attempts were presented in these tables, with no malicious intent.

11. *« And what is the point of additions where, for example, amalgamating three sets of 3, 6 and 4 into a single set of 13 does nothing but make the number more difficult to understand? »*

Only that it brings a visual illustration of the sum of the three sets of 3, 6 and 4 notches on the middle column as a new set of 13 notches in the left column.

**Intention of prehistoric men and unknown**

12. « *It is well known that such an "addition" would have been completely meaningless in the first true number systems. To return to the case in hand, the number 13 would never have been represented with 13 regularly spaced notches, but instead in distinct sets to make it easier to understand. In Ancient Egypt, for example, the hieroglyph 9 was not 9 equally spaced aligned bars, but either 4 bars placed below 5 other bars, or more often three sets of 3 bars placed atop one another.* »

Is it a criticism or an argument? The authors no longer understand this critic, because on one hand, he decries their similar interpretations, and on the other hand, he proposes himself the bold idea of the passage of a representation of a set of 13 notches to a more abstract representation of a number, 13. And that was certainly not the intention of prehistoric men of Ishango, who had probably not yet discovered the concept of numbers and their representation. This passage to abstraction and representation of the numbers is well posterior to the prehistoric civilization of Ishango.

In addition, this critic points out in his footnote 9 of the page 9 of [19], the Babylonian numeration using a symbol (the nail) in base 10 and another (the chevron) in base 6. The number 9 is then written with nine nails, according to B. Rittaud [26, p.2], quoted by this critic. So what: two weights, two measures?

13. « *It doesn't work, so let's invoke unknown intentions! These are joined by some rather unconvincing speculations about the comparative length or gradient of the notches, which in any case do not justify the necessary addition of 2, 1 or -1 to lend a semblance of coherence to the whole.* »

Regarding the « *unknown intentions* » in the first sentence, the authors signal that it happens also to this critic himself to confess his ignorance, for example [27, p. 73]:

> « *For the reasons that the author ignores, the primitive symbols are essentially geometric symbols.* »

In addition, these « *unconvincing speculations* » refer to the study of the shape and size of the notches in the various columns that are detailed in [15].

**The choice of numbers on the bone**

14. « *As for the fifth number in the middle column, one might well ask why the authors choose 10, which does not give a single correct answer, rather than 9, which gives four correct answers. The reason is that with 10, the total of the middle column is 48, which is a multiple of 12, like the total 60 in the left and right columns.* »

Indeed.

15. « *To account for the choice of numbers on the bone, Pletser and Huylebrouck posit that: « The numbers 3 and 4 could have formed the base of the arithmetic system used by the ancient Ishango people for operations on small numbers and that the derived base 12 was used for larger numbers. » Where are the bases 3 and 4? In the middle column, according to the authors, because from top to bottom it shows:*

*- 3 then 6, so 3 then 3 x 2*
*- 4 then 8, so 4 then 4 x 2*
*- 9 or 10, so 4 x 2+1 or 4 x 2+2*
*- two times 5, « showing two ways of obtaining the composed number 5, in adding 1 or 2 to either of the bases 3 and* [sic] *4 »*
*- 7 « showing how to obtain the composed number 7 by adding the two bases 3 and 4»*

Indeed.

**Combination of several bases**

16. *« But if 5 is shown twice because there are two bases, why are the other numbers – 6, 8, 9 or 10 and 7 – shown only once? Furthermore, if the aim had been to demonstrate a base, this would have been clear to see. Two groups of 3 should be visible within a set of 6, two groups of 4 should be visible within a set of 8, and so on and so forth. Yet there is nothing of the kind. No regular groupings can be detected that are suggestive of a base. »*

Nowhere was it said that the sculptor of notches had *«the aim to demonstrate a base»* on this bone. The authors do not understand where this critic's remark comes from. The authors suggested that this rod is probably the testimony of a people who were counting by mixing bases 10 and 12 (or 6). A European would maybe write ||||| ||||| |||| to count 14 days, whereas a man from Ishango may write |||||| ||||||| ||. The combination of several bases is common and is not surprising. For example, in French, 77 and 93 refer to a past in bases 10 and 20. If in 22,000 years, one will find a French text with the expressions 'seventy' and 'eighty-three[2]', the archaeologist from the future might conclude with reason that the French language was mixing bases 10 and 20. In addition, there would be little chance that this text would be a scholarly dissertation on arithmetic explaining the use of these different bases.

Then, the duplication is not necessarily an operation of addition to the original of a copy of the original, or a simple multiplication by 2. It may be a reconstruction of a new group with a number of notches double of the initial group but arranged differently [15]. Requiring *« to see clearly two sets of 3 within the set of 6, two sets of 4 within the set of 8 and so forth »* is what might be called 'looking through the small end of the spyglass' and demonstrates a too simplistic modern mathematical-cultural egocentricity.

17. *« Where is the base 12? On one hand, as we have already noted, in the column totals, which are multiples of 12. And on the other hand, according to the authors, in the fact that in the middle column:*
*- 6 is involved in two assumed additions (first two lines of the above table): 2 x 6 = 12*
*- 4 is involved in three assumed additions: 3 x 4 = 12*
*- 8 is involved in three assumed additions: 3 x 8 = 24 = 2 x 12»*

Indeed, but, as mentioned above, the base 12 also appears in the left column with the numbers 11, 13, 17, 19, which, after the abandonment of the hypothesis of the prime numbers, are seen as $12 \pm 1$ and $18 \pm 1$.

18. *« This gives us the following situation: 6 is incised once as a group of six notches in the middle column. But, because we have assumed that it is involved in two additions, and even*

---

[2] In French, the number 83 is said and written literally 'four-twenty-three' ('*quatre-vingt-trois*').

*though these are incorrect (vide the unknown intention), that gives 12! The same goes for 4 and 8, which supposedly appear three times each, giving 12 and 2 x 12 respectively. How can one possibly be convinced by such sleights of hand as these, which are worthy of the most insipid numerological tract? »*

Let's put things in their place. The bases 3 and 4 are deduced from the basic operations in the middle column and the base 12 is deduced from the numbers 12 ± 1 and 18 ± 1 of the left column. Simply, the authors find that the numbers of the middle column are equal to 12 or multiples or submultiples of 12. This fact is nowhere used as an argument to reinforce the veracity of the hypothesis of the base 12 and its submultiples.

**The second bone**

19. *« In 1959, Jean de Heinzelin found another notched bone, again in Ishango. In 1998 he put forward an interpretation of these notches which, according to Pletser and Huylebrouck, confirms the above. However, there is little point in continuing to test the reader's patience with this matter. »*

The authors doubt that this critic has seen the second rod, because nothing had been published on this subject prior to 2007, on the occasion of the Congress to which this critic was invited, but that he preferred to avoid.
The authors believe instead that this point is of great interest. On his deathbed, Jean de Heinzelin had changed his mind about the interpretation of the notches of the first rod, actually seeing a basic arithmetic tool using several simple bases. This conclusion had not been published according to his last wishes, and would be published only some time after his death. The authors were very surprised to see their hypothesis confirmed *a posteriori* by de Heinzelin, at the same time as the announcement of the existence of the second bone.

20. *« A glance at Figure 5 and the following passage will be edifying enough:*

> *« Prof. De Heinzelin added that the minor on the E Column is at the "10- spot", and wondered if this announced "a passage from the base 10 to base 12" […] Since the C column has a total of 20 carvings, and the E column 18 […] the bases 6 and 10–20 seem to emerge. Moreover, there are two spatial concordances between the rows, at E10 = F1 = G10 and at E12 = F2 = G12.12. »*

*According to de Heinzelin's hypotheses, which were later taken up by Pletser and Huylebrouck, this bone may have borne witness to a change of base, or had a didactic function, or even played a role in exchanges between different ethnic groups, some of which used base 10, while others used bases 12 or 16 among others. »*

The archaeologist de Heinzelin was a world leader in archaeology. His arguments were always based and had weight.

**The small end of the spyglass**

21. *« But remember, the base argument can be taken seriously only if the groupings are clear and systematic. 18 = 3 x 6 is not proof of base 6! »*

And remember also that this way of seeing things is contrary to the scientific approach which presupposes an open mind and the examination of all possible and plausible hypotheses.

Looking only for « *clear and systematic groupings* » is so reducing and simplistic, and is evidence of a narrowness of view, (the small end of the spyglass...) passing any observation through the mould of our modern vision of things without questioning the way in which a prehistoric man may have seen things.

Indeed, $18 = 3 \times 6$ is not enough to deduce the existence of a base 6 if there is only one indication of this type. On the other hand, if several indications are present on the first rod of Ishango as shown by de Heinzelin and other scientists, with several examples of basic operations of duplication and addition involving the numbers 3, 4 and 12, next to the number 10, any scientist who has been trained in the scientific method, namely observation and deduction, is entitled to ask the question if there is actually more than a simple chance and if correlations between the different bases should not be considered.

**Ethnographic comparatism**

22. *«The handful of ethnographic examples adduced by the authors of the conference proceedings are just as unconvincing. The fact that the people of the Congo say the equivalent of "twelve-one" when they mean thirteen makes base 12 relevant here, but what it actually signifies is that the only way to say 13 is 12, then 1. The visual equivalent would be to incise a set of 12 marks followed by a space and a separate notch. The same goes for peoples who use their fingers to represent numbers. The Shambaa of Tanzania represent the number 6 by stretching out three fingers on each hand, and say the equivalent of "three- three" for six. Number 8 is "four-four" and represented by four fingers on each hand. Number 7 is more complex, in that it is pronounced as "ten minus three" and represented by four fingers on the right hand and three on the left hand. The gestures for the three numbers – 7, 8 and 6 – make a clear distinction between bases 3 and 4, if the term base is indeed appropriate here. Yet such a clear and systematic separation into subsets of 3, 4 and 12 notches is not in evidence on either of the Ishango bones. These ethnographic examples only make matters worse for Pletser and Huylebrouck's theories. »*

It is surprising to read the beginning of this excerpt from the pen of someone who presents himself as the promotor of the ethnographic comparatism, and who is even ecstatic [22, pp. 563-564]:

> « *The ethnographic comparatism, which recognizes this analogy, opens a research path barely explored for the mathematics of prehistory, but rich of promises since it allows, to a certain extent, to make talk the archaeological findings. And it also imposes on the assumptions and theories based on the prehistoric documents to be verified by the ethnographic documents and vice versa.* »

The authors do not see how or why « *These ethnographic examples only make matters worse for the theories* » of the authors. There is nothing that allows to say this, as this critic does. On the contrary, if he had bothered to read carefully the different ethnographic circumstantial arguments given in [15].

**A question of battered women?**

23. « *Yet, whether simple or sophisticated, such purported interpretations all have the same arbitrary foundation: the belief that notches are necessarily numerical. Claudia Zaslavsky recounts that some African women occasionally make a notch on the handle of their wooden spoon. Are they marking the passing days? Or playing with numbers? Not all* [sic]*: they make* *a notch each time their husband hits them, and when the spoon handle is full, they ask for a divorce.* »

The authors cited Zaslavsky in their reference 23, p. 166 [15]. More even, they cited similar reports written by missionaries and administrators of the Belgian Congo who described how and why such rods with notches were used.

In addition, it is true that this hypothesis may not be entirely rejected. But then, why an Ishango prehistoric woman would have made marks so sophisticated with groups that provide such regularity in the whole? And moreover, why make marks in three columns? One can imagine that she could have already asked for a divorce after filling the middle column with 48 notches, or the left or the right column, with 60 notches each time!

Although this hypothesis cannot be formally rejected, it is clear that this hypothesis would be the least likely of all for the first rod of Ishango.

24. « *A notch may be nothing more than a mark, which seems like small fry if one is obsessed with arithmetic. And yet that is the most important invention we owe to our ancestors of the Upper Palaeolithic* [sic]*: the sign* »

Indeed.

25. « *By losing ourselves in haphazard mathematical speculations, we waste time, money and paper, and this when there is so much to discover about prehistoric signs – including the intellectual gestation of the concept of number – by considering them alongside the ethnographic records.* »

If this critic had followed this recommendation to the letter, his scientific work, confined to a few publications at the national level and repeating the same criticisms as those discussed here, would have disappeared entirely.

26. The two ethnographic counterexamples given by this critic [19, pp.14 - 17] are excellent but are not counter-examples, on the contrary. The authors take similar examples in some of their publications. These stories show especially an important point: the *a priori* knowledge of the reason for being, or *raison d'être*, of these marked or signed media. In other words, another source of information (written, oral tradition, additional information, whatever the medium) allows an easy decryption of these objects.

The Ishango bones have been found with other tools and spear tips [21]. None of these other artefacts indicate any link with the notches on the bones. The critic says himself [22, p. 567] speaking of the Ishango bone:

> « *Similarly, no ethnographic document allows to support the above-mentioned thesis [...] of a hunter-gatherer producing a table of primes, a table of doubles, ...* »

So what to do? Ignore these artefacts and leave them in a drawer or on the contrary study them and try to decipher them?

Remember that the study of the first bone was done without even knowing that a second bone had been discovered in the same place and that the description of this second bone has been published several years after the death of its discoverer. Yet the conclusions drawn on the first bone alone by the authors have been corroborated afterward by the conclusions of Jean de Heinzelin before the publication of the authors' results and without knowledge of de Heinzelin conclusions.

**Which calculator?**

27. « *The currency of such fictions as the "Ishango calculator" – and the fact that they are often taken at face value – is a sorry state of affairs, not only because of their intrinsic flimsiness and implausibility, but also because the archaeological and ethnographic archives could be put to so much better use.* »

The authors have never used the term *«Ishango calculator»* in their publications and communications. They leave the responsibility and paternity of this expression to this critic.

## Conclusions:

**An appropriate conclusion**

The study of the first bone led to an appropriate conclusion, although not spectacular. It was reached without knowing that a second bone had been discovered in the same place and before the publication of its description several years after the death of its discoverer. Yet the conclusions that the authors drew on the first bone have been corroborated *a posteriori* by Jean de Heinzelin's conclusions on his deathbed before the publication of the results of the authors and without knowledge of de Heinzelin's conclusions.

**An epistemological fault**

A better bibliographic research would have allowed this critic to forge his opinion on all the elements, but one can doubt that this critic has a sufficient command of the English language, as he pulls items out of their context and uses second hand sources instead of the original documents. Well, it is true that this is not the first time that this critic has problems to read scientific publications and to document himself correctly (see for example [28]).

**Breach of copyright**

The authors also signal the use by this critic of images without respect for copyright laws; his illustrations are taken from other articles without any authorization request, or from the site of the Museum of Natural Sciences in Brussels, without respect for applicable regulations.

**A doubtful, little academic vocabulary**

The authors can't resist the pleasure to highlight in this little pamphlet the adjectives and expressions which shed light on the style of this critic: "*ridiculous*" ([19], p. 5), "*they have to forcibly make up others*" (p. 6), "*a muddled table*" (p. 7), "*It doesn't work*" (p. 8), "*unconvincing speculations*" (p. 8), "*a semblance of coherence*" (p. 8), "*which are worthy of the most insipid numerological tract?*" (p. 9), "*this "trick" doesn't work*" (p. 9), "*Since we need to have 12s*" (p. 9), "*such cobbled-together computations*" (p. 9), "*to test the reader's patience with this matter. A glance at Figure 5 and the following passage will be edifying enough*" (p. 10), "*Next comes the invention of some kind of concrete context and a story. In the end, what we have before us is mathematical fiction.*" (p. 11), "*speculations*" (p. 11), "*the latest in a long line of victims of* *mathematical illusion*" (p.11), "*The siren song of mathematical illusion*" (p. 11), "*a naïve mathematician*" (p. 11), "*Four is feminine because women have four lips*." (Note 24, p. 14), "*calamitous enough for mathematical fablers*" (p. 15).

As for the terms "*intrinsic flimsiness and implausibility*" (p. 16), they would be better attached to this little critical pamphlet.

The authors deplore this language but they do not want to get into a pointless controversy for reasons explained in the introduction.

However, they are open to debating ideas and they would be happy to participate in an exchange of views at a conference.

*V. Pletser, D. Huylebrouck, August 2015*

References

[1] D. Huylebrouck, 'Puzzles, patterns, drums: the dawn of mathematics in Rwanda and Burundi', *Humanistic Mathematics Network, Journal* no 14, November 1996.

[2] D. Huylebrouck, 'The Bone that began the Space Odyssey', *The Mathematical Intelligencer*, Vol. 18, no 4, p. 56, Autumn 1996.

[3] D. Huylebrouck, 'Histoires célestes du Rwanda' (Celestial stories in Rwanda), *Ciel et Terre, Bulletin of the Royal Belgian Soc. of Astronomy, Meteorology and Geophysics*, Vol. 113 (2), 83-89, 1997.

[4] D. Huylebrouck, 'Counting on hands in Africa and the origin of the duodecimal system', *Wiskunde en Onderwijs (Mathematics and Education)*, no 89, 1997.

[5] V. Pletser, D. Huylebrouck, 'The Ishango Artefact: the Missing Base 12 Link.' *Proc. Katachi Univ. Symmetry Congress (KUS2)*, T. Ogawa, S. Mitamura, D. Nagy & R. Takaki (ed.), Paper C11, Tsukuba Univ., Japan, 18 Nov. 1999.

[6] V. Pletser, D. Huylebrouck, 'Research and promotion about the first mathematical artefact: the Ishango bone', *Proc. PACOM 2000 Meeting Ethnomathematics and History of Mathematics in Africa*, Cape Town, South Africa, 1999.

[7] V. Pletser, D. Huylebrouck, 'The Ishango Artefact: the Missing Base 12 Link', *Forma*, 14-4, 339-346, 1999.

[8] D. Huylebrouck, '*Afrika + Wiskunde*' (in Dutch), VUBPress of the Free University of Brussels, 2005; translated in French '*Afrique + Mathématiques*' VUBPress of the Free University of Brussels, 2008.

[9] D. Huylebrouck, 'De parabolische vlucht van het Ishangoartefact, de oudste wiskundige vondst' (The parabolic flight of the Ishango artefact, the oldest mathematical discovery), *Wiskunde en Onderwijs (Mathematics and Education)*, n° 101, 2000.

[10] D. Huylebrouck, F. Dumortier, 'A Stamp for World Mathematical Year', *European Mathematical Society Newsletter*, p. 27, November 2000.

[11] D. Huylebrouck, 'Middle column of marks found on the oldest object with logical carvings, the 22000-year-old Ishango bone from the Congo', sequence A100000, *On-Line Encyclopaedia of Integer Sequences*, http://oeis.org/, November 2004

[12] D. Huylebrouck, 'L'Afrique, berceau des mathématiques' (Africa, the craddle of mathematics), *Pour la Science*, Dossier 47, April-June 2005 (6 pages). The article was translated in Portuguese (in Brazil) and in German.

[13] D. Huylebrouck, 'Mathematics in (central) Africa before colonisation', *Anthropologica et Praehistorica, Bulletin of the Royal Belgian Association for Anthropology and Prehistory*, 135 – 162, Vol. 117, 2006.

[14] D. Huylebrouck, 'Report: The ISShango project', *Journal for Mathematics and the Arts*, Taylor and Francis, September 2008.

[15] V. Pletser, D. Huylebrouck, 'An interpretation of the Ishango rods', *Proc. Conf. 'Ishango, 22000 and 50 years later: the cradle of Mathematics?'*, D. Huylebrouck ed., Royal Flemish Academy of Belgium, KVAB, 139-170, 2008.

[16] D. Huylebrouck, '*Foreword to the book Lusona: find the omitted drawings*' by Paulus Gerdes (Mozambique), Lulu Editions, July 2009.

[17] D. Huylebrouck, R. Matteus Berr, 'Vermessung', in '*Kunstlerhaus catalogue EVO EVO! 200 Jahre Charles Darwin*', I. and P. Braunsteiner eds, Künstlerhaus ISBN 978-3-900926-84-7 Verlag Lehner Wien, 2009.

[18] D. Huylebrouck, '*België + wiskunde' (Belgique + Mathématique*), Academia Press Gent, Belgique, 2013.

[19] O. Keller, « The fables of Ishango, or the irresistible temptation of mathematical fiction », *Bibnum*, August 2010.

[20] J. Yellen, 'Ishango and its importance for later African prehistory', *Proc. Conf. 'Ishango, 22000 and 50 years later: the cradle of Mathematics?'*, D. Huylebrouck ed., Royal Flemish Academy of Belgium, KVAB, 63-81, 2008.

[21] J. de Heinzelin de Braucourt, 'Exploration du Parc National Albert', Fascicule 2, *Les Fouilles d'Ishango*, Instituts des Parc Nationaux du Congo Belge, Bruxelles, 1957.

[22] O. Keller, « Questions Ethnographiques et Mathématiques de la Préhistoire », *Revue de synthèse* : 4e S. n° 4, oct.-dec. 1998, p. 545-573.

[23] O. Keller, « Préhistoire de l'arithmétique, La découverte du nombre et du calcul», in *'Si le nombre m'était conté'*, I.R.E.M. - Histoire des mathématiques, Paris, Ellipse, 2000.

[24] F.P. Ramsey, 'On a Problem of Formal Logic', *Proceedings London Mathematical Society*, s2-30 (1): 1930, 264–286. doi10.1112/plms/s2-30.1.264

[25] L. Hogben, «*Mathematics for the Million*», London, Allen & Unwin, 1936. ISBN-13: 978-0393310719

[26] B. Rittaud, « À un mathématicien inconnu ! », *Bibnum*, novembre 2008.

[27] O. Keller, « Préhistoire de la Géométrie: Premiers éléments d'enquête, premières conclusions », *Science et Techniques en Perspective*, Vol. 33, Université de Nantes, 1995.

[28] J. Høyrup, 'Book review of 'Préhistoire de la Géométrie : Premiers éléments d'enquête, premières conclusions', published in *Isis* 87, 1996, 713-714. Available on http://www.academia.edu/3131829/Prehistoire_de_la_Geometrie_Premiers_elements_denquete_premieres_conclusions
http://akira.ruc.dk/~jensh/publications/1996%7BR%7D17_Keller_Prehistoire_MS.PDF